\documentclass[letterpaper,USenglish]{article}

\usepackage{amssymb}
\usepackage{graphicx}
\usepackage{slashbox}
\usepackage{pict2e}
\usepackage{hyperref}
\usepackage{caption}
\newtheorem{theorem}{Theorem}

\newtheorem{lemma}{Lemma}
\newtheorem{definition}{Definition}

\def\pont{\hspace{-6pt}{\bf.\ }}
\def\eps{\varepsilon}
\title{A Practical Regularity Partitioning Algorithm and its Applications in Clustering}

\author{G\'{a}bor N. S\'{a}rk\"{o}zy\thanks{Research supported in part by
by the National Science Foundation under Grant No. DMS-0968699.} \dag, Fei Song\dag, Endre Szemer\'edi\ddag, Shubhendu Trivedi\dag \\ \\
\dag Computer Science Department\\Worcester
Polytechnic Institute\\Worcester, MA 01609 \\ \\
\ddag Computer Science Department \\ Rutgers University \\
New Brunswick, NJ 08901}

\begin{document}

\maketitle

\begin{abstract}
In this paper we introduce a new clustering technique called {\em
Regularity Clustering}. This new technique is based on the practical
variants of the two constructive versions of the Regularity Lemma, a
very useful tool in graph theory. The lemma claims that every graph
can be partitioned into pseudo-random graphs. While the Regularity
Lemma has become very important in proving theoretical results, it
has no direct practical applications so far. An important reason for
this lack of practical applications is that the graph under
consideration has to be astronomically large. This requirement makes
its application restrictive in practice where graphs typically are
much smaller. In this paper we propose modifications of the
constructive versions of the Regularity Lemma that work for smaller
graphs as well. We call this the Practical Regularity partitioning
algorithm. The partition obtained by this is used to build the
reduced graph which can be viewed as a compressed representation of
the original graph. Then we apply a pairwise clustering method such
as spectral clustering on this reduced graph to get a clustering of
the original graph that we call Regularity Clustering. We present
results of using Regularity Clustering on a number of benchmark
datasets and compare them with standard clustering techniques, such
as $k$-means and spectral clustering. These empirical results are
very encouraging. Thus in this paper we report an attempt to harness
the power of the Regularity Lemma for real-world applications.
 \end{abstract}

\section{Introduction}

The Regularity lemma of Szemer\'edi \cite{Sz} has proved to be a
very useful tool in graph theory. It was initially developed as an
auxiliary lemma to prove a long standing conjecture of Erd\H{o}s and
Tur\'{a}n\cite{ET} on arithmetic progressions, which stated that
sequences of integers with positive upper density must contain
arbitrarily long arithmetic progressions. Now the Regularity Lemma
by itself has become an important tool and found numerous other
applications (see \cite{KSSS}). Based on the Regularity Lemma and
the Blow-up Lemma \cite{KSS5} the Regularity method has been
developed that has been quite successful in a number of applications
in graph theory (e.g. \cite{GRSS1}). However, one major disadvantage
of these applications and the Regularity Lemma is that they are
mainly theoretical, they work only for astronomically large graphs
as the Regularity Lemma can be applied only for such large graphs.
Indeed, to find the $\varepsilon$-regular partition in the
Regularity Lemma, the number of vertices must be a tower of 2's with
height proportional to $\varepsilon^{-5}$. Furthermore, Gowers
demonstrated \cite{G} that a tower bound is necessary.

The basic content of the Regularity Lemma could be described by
saying that every graph can, in some sense, be partitioned into
random graphs. Since random graphs of a given edge density are much
easier to treat than all graphs of the same edge-density, the
Regularity Lemma helps us to carry over results that are trivial for
random graphs to the class of all graphs with a given number of
edges. We are especially interested in harnessing the power of the
Regularity Lemma for clustering data. Graph partitioning methods for
clustering and segmentation have become quite popular in the past
decade because of representative ease of data with graphs and the
strong theoretical underpinnings that accompany the same.

In this paper we propose a general methodology to make the
Regularity Lemma more useful in practice. To make it truly
applicable, instead of constructing a provably regular partition we
construct an {\em approximately} regular partition. This partition
behaves just like a regular partition (especially for graphs
appearing in practice) and yet it does not require the large number
of vertices as mandated by the original Regularity Lemma. Then this
approximately regular partition is used for performing clustering.
We call the resulting new clustering technique {\em Regularity
clustering}. We present comparisons with standard clustering methods
such as $k$-means and spectral clustering and the results are very
encouraging.

To present our attempt and the results obtained, the paper is
organized as follows: In section \ref{prior} we discuss briefly some
prior attempts to apply the Regularity Lemma in practical settings
and place our work in contrast to those. In Section \ref{clust} we
discuss clustering in general and also present a popular spectral
clustering algorithm that is used later on the reduced graph. We
also point out what are the possible ways to improve its running
time. In Section \ref{not} we give some definitions and general
notation. In Section \ref{reg} we present two constructive versions
of the Regularity Lemma (the original lemma was non-constructive).
Furthermore, in this section we point out the various problems
arising when we attempt to apply the lemma in real-world
applications. In Section \ref{mod} we discuss how the constructive
Regularity Lemmas could be modified to make them truly applicable
for real-world problems where the graphs typically are much smaller,
say have a few thousand vertices only. In Section \ref{app} we show
how this Practical Regularity partitioning algorithm can be applied
to develop a new clustering technique. In Section \ref{test}, we
present an extensive empirical validation of our method. Section
\ref{future} is spent in discussing the various possible future
directions of work.

\section{Prior Applications of the Regularity Lemma}\label{prior}

As we discussed above so far the Regularity Lemma has been ``well
beyond the realms of any practical applications" \cite{Abel}, the
existing applications have been theoretical, {\em mathematical}. The
only practical application attempt of the Regularity Lemma to the
best of our knowledge is by Sperotto and Pelillo \cite{SP}, where
they use the Regularity Lemma as a pre-processing step. They give
some interesting ideas on how the Regularity Lemma might be used,
however they do not give too many details. Taking leads from some of
their ideas we give a much more thorough analysis of the
modifications needed in order to make the lemma applicable in
practice. Furthermore, they only give results for using the
constructive version by Alon {\em et al} \cite{AD}, here we
implement the version proposed by Frieze and Kannan \cite{FK} as
well. We also give a far more extensive empirical validation; we use
12 datasets instead of 3.

\section{Clustering}\label{clust}

Out of the various modern clustering techniques, spectral clustering
has become one of the most popular. This has happened due to not
only its superior performance over the traditional clustering
techniques, but also due to the strong theoretical underpinnings in
spectral graph theory and its ease of implementation. It has many
advantages over the more traditional clustering methods such as
$k$-means and expectation maximization (EM). The most important is
its ability to handle datasets that have arbitrary shaped clusters.
Methods such as $k$-means and EM are based on estimating explicit
models of the data. Such methods fail spectacularly when the data is
organized in very irregular and complex clusters. Spectral
clustering on the other hand does not work by estimating explicit
models of the data but does so by analysing the spectrum of the
Graph Laplacian. This is useful as the top few eigenvectors can
unfold the data manifold to form meaningful clusters.

In this work we employ spectral clustering on the reduced graph (which is an essence
of the original graph),
even though any other pairwise clustering method could be used. The
algorithm that we employ is due to Ng, Jordan and Weiss \cite{NJW}.
Despite various advantages of spectral clustering, one major problem
is that for large datasets it is very computationally intensive. And
understandably this has received a lot of attention recently. As originally stated,
the spectral clustering pipeline has two main bottlenecks: First, computing the
affinity matrix of the pairwise distances between datapoints, and second, once we have the
affinity matrix the finding of the eigendecomposition. Many ways have been
suggested to solve these problems more efficiently. One approach is not to
use an all-connected graph but a k-nearest neighbour graph in which
each data point is typically connected to $\log{n}$ neighboring
datapoints(where $n$ is the number of data-points). This
considerably speeds up the process of finding the affinity matrix,
however it has a drawback that by taking nearest neighbors we might
miss something interesting in the global structure of the data. A
method to remedy this is the Nystr\"{o}m method which takes a random
sample of the entire dataset (thus preserving the global structure
in a sense) and then doing spectral clustering on this much smaller
sample. The results are then extended to all other points in the
data set \cite{FowlChung}.

Our work is quite different from such methods. The speed-up is
primarily in the second stage where eigendecomposition is to be
done. The original graph is represented by a reduced graph which is
much smaller and hence eigendecomposition of this reduced graph can significantly
ease the computational load. Further work on a practical variant
of the sparse Regularity Lemma could be useful in a speed-up in the
first stage, too.

\section{Notation and Definitions}\label{not}

Below we introduce some notation and definitions for describing the
Regularity Lemma and our methodology.

Let $G = (V,E)$ denote a graph, where $V$ is the set of vertices and
$E$ is the set of edges. When $A, B$ are disjoint subsets of $V$,
the number of edges with one endpoint in $A$ and the other in $B$ is
denoted by $e(A,B)$. When $A$ and $B$ are nonempty, we define the
{\em density} of edges between $A$ and $B$ as $d(A,B) =
\frac{e(A,B)}{|A||B|}$. The most important concept is the following.

\begin{definition}\pont
The bipartite graph $G=(A,B,E)$ is $\varepsilon$-{\em regular} if
for every $X\subset A$, $Y\subset B$ satisfying: $\
|X|>\varepsilon|A|,\ |Y|>\varepsilon|B|,$ we have
$|d(X,Y)-d(A,B)|<\varepsilon,$ otherwise it is $\varepsilon$-{\em
irregular}.
\end{definition}

\noindent Roughly speaking this means that in an $\varepsilon$-regular bipartite
graph the edge density between {\em any} two relatively large
subsets is about the same as the original edge density. In effect
this implies that all the edges are distributed almost uniformly.

\begin{definition}\pont\label{defn}
A partition $P$ of the vertex set $V=V_0\cup V_1\cup \ldots \cup
V_k$ of a graph $G = (V,E)$ is called an {\em equitable partition}
if all the classes $V_i, 1\leq i\leq k$, have the same cardinality.
$V_0$ is called the exceptional class.
\end{definition}

\begin{definition}\pont\label{potential}
For an equitable partition $P$ of the vertex set $V=V_0\cup V_1\cup
\ldots \cup V_k$ of $G = (V,E)$, we associate a measure called the
{\em index} of $P$ (or the potential) which is defined by
$$ind(P) = \frac{1}{k^2}\sum_{s=1}^k \sum_{t=s+1}^k d(C_s,C_t)^2.$$
\end{definition}
This will measure the progress towards an $\varepsilon$-regular partition.

\begin{definition}\pont
An equitable partition $P$ of the vertex set $V=V_0\cup V_1\cup
\ldots \cup V_k$ of $G = (V,E)$ is called $\varepsilon$-{\em
regular} if $|V_0| < \varepsilon |V| $ and all but $\varepsilon k^2$
of the pairs $(V_i,V_j)$ are $\varepsilon$-regular where $1 \leq i <
j \leq k$.
\end{definition}
With these definitions we are now in a position to state the
Regularity Lemma.

\section{The Regularity Lemma}\label{reg}

\begin{theorem}[Regularity Lemma \cite{Sz}]\pont\label{p1}
For every positive $\varepsilon > 0$ and positive integer $t$ there is an
integer $T = T(\varepsilon,t)$ such that every graph with $n > T$ vertices
has an $\varepsilon$-regular partition into $k+1$ classes, where $t \leq k
\leq T$.
\end{theorem}

In applications of the Regularity Lemma the concept of the {\em
reduced graph} plays an important role.

\begin{definition}\pont\label{reduced}
Given an $\varepsilon$-regular partition of a graph $G = (V, E)$ as
provided by Theorem \ref{p1}, we define the {\em reduced graph}
$G^R$ as follows. The vertices of $G^R$ are associated to the
classes in the partition and the edges are associated to the
$\varepsilon$-regular pairs between classes with density above $d$.
\end{definition}

The most important property of the reduced graph is that many
properties of $G$ are inherited by $G^R$. Thus $G^R$ can be treated
as a representation of the original graph $G$ albeit with a much
smaller size, an ``essence'' of $G$. Then if we run any algorithm on
$G^R$ instead of $G$ we get a significant speed-up.

\subsection{Algorithmic Versions of the Regularity Lemma} \label{RegAlgo}

The original proof of the Regularity Lemma \cite{Sz} does not give a
method to construct a regular partition but only shows that one must
exist. To apply the Regularity Lemma in practical settings, we need
a constructive version. Alon {\em et al.} \cite{AD} were the first
to give an algorithmic version. Since then a few other algorithmic
versions have also been proposed \cite{FK}, \cite{KRT}. Below we
present the details of the Alon {\em et al.} algorithm.

\subsubsection{Alon {\em et al.} Version}

\begin{theorem}[Algorithmic Regularity Lemma \cite{AD}]\pont\label{al}
For every $\varepsilon > 0$ and every positive integer $t$ there is an
integer $T = T(\varepsilon, t)$ such that every graph with $n > T$ vertices
has an $\varepsilon$-regular partition into $k + 1$ classes, where $t \le k
\le T$. For every fixed $\varepsilon > 0$ and $t \ge 1$ such a partition
can be found in $O(M(n))$ sequential time, where $M(n)$ is the time
for multiplying two $n$ by $n$ matrices with $0, 1$ entries over the
integers. The algorithm can be parallelized and implemented in
$NC^1$.
\end{theorem}

This result is somewhat surprising from a computational complexity
point of view since as it was proved in \cite{AD} that the
corresponding decision problem (checking whether a given partition
is $\varepsilon$-regular) is co-NP-complete. Thus the search problem
is easier than the decision problem. To describe this algorithm, we
need a couple of lemmas.

\begin{lemma}[Alon {\em et al.} \cite{AD}]\pont\label{lnew1}
Let $H$ be a bipartite graph with equally sized classes $|A| = |B| =
n$. Let $2n^{-1/4} < \varepsilon <\frac{1}{16}$. There is an $O(M(n))$
algorithm that verifies that $H$ is $\varepsilon$-regular or finds two
subset $A' \subset A$, $B' \subset B$, $|A'| \ge
\frac{{\varepsilon}^4}{16}n$, $|B'| \ge \frac{{\varepsilon}^4}{16}n$, such that
$|d(A, B) - d(A', B')| \ge \varepsilon^4$. The algorithm can be
parallelized and implemented in $NC^1$.
\end{lemma}

This lemma basically says that we can either verify that the pair is
$\varepsilon$-regular or we provide certificates that it is not. The
certificates are the subsets $A', B'$ and they help to proceed to
the next step in the algorithm. The next lemma describes the
procedure to do the refinement from these certificates.

\begin{lemma}[Szemer\'edi \cite{Sz}]\pont\label{lnew2}
Let $G = (V,E)$ be a graph with $n$ vertices. Let $P$ be an
equitable partition of the vertex set $V=V_0\cup V_1\cup \ldots \cup
V_k$. Let $\gamma >0$ and let $k$ be a positive integer such that
$4^k > 600\gamma^{-5}$. If more than $\gamma k^2$ pairs $(V_s,
V_t)$, $1 \le s < t \le k$, are $\gamma$-irregular then there is an
equitable partition $Q$ of $V$ into $1 + k4^k$ classes, with the
cardinality of the exceptional class being at most $|V_0| +
\frac{n}{4^k}$ and such that $ind(Q) > ind(P) + \frac{\gamma^5}{20}.$
\end{lemma}

This lemma implies that whenever we have a partition that is not
$\gamma$-regular, we can refine it into a new partition which has a
better index (or potential) than the previous partition. The
refinement procedure to do this is described below.

{\bf Refinement Algorithm:} {\em Given a $\gamma$-irregular
equitable partition $P$ of the vertex set $V=V_0\cup V_1\cup \ldots
\cup V_k$ with $\gamma = \frac{\varepsilon^4}{16}$, construct a new partition $Q$.\\
For each pair $(V_s, V_t)$, $1 \leq s < t \leq k$, we apply Lemma
\ref{lnew1} with $A=V_s$, $B=V_t$ and $\varepsilon$. If $(V_s, V_t)$
is found to be $\varepsilon$-regular we do nothing. Otherwise, the
certificates partition $V_s$ and $V_t$ into two parts (namely the
certificate and the complement). For a fixed $s$ we do this for all
$t\not= s$. In $V_s$, these sets define the obvious equivalence
relation with at most $2^{k-1}$ classes, namely two elements are
equivalent if they lie in the same partition part for every $t\not=
s$. The equivalence classes will be called atoms. Set $m = \lfloor
\frac{|V_i|}{4^k}\rfloor$, $1 \le i \le k$. Then we construct our
new partition $Q$ by choosing a maximal collection of pairwise
disjoint subsets of $V$ such that every subset has cardinality $m$
and every atom $A$ contains exactly $\lfloor \frac{|A|}{m}\rfloor$
subsets; all other vertices are put in the exceptional class. The
collection $Q$ is an equitable partition of $V$ into at most
$1+k4^k$ classes and the cardinality of its exceptional class is at
most $|V_0| + \frac{n}{4^k}$. }

Now we are ready to present the main algorithm.

{\bf Regular Partitioning Algorithm:} {\em Given a graph $G$ and
$\varepsilon$, construct a $\varepsilon$-regular partition.
\begin{enumerate}
\item {\bf Initial partition:} Arbitrarily divide the vertices of $G$ into an equitable partition $P_1$ with classes $V_0, V_1, \ldots, V_b$, where
$|V_1| = \lfloor \frac{n}{b} \rfloor$ and hence $|V_0| < b$. Denote
$k_1 = b$.
\item {\bf Check regularity:} For every pair $(V_s,V_t)$ of $P_i$, verify if it is $\varepsilon$-regular or find
$X \subset V_s, Y \subset V_t, |X| \ge \frac{\varepsilon^4}{16}|V_s|,|Y|
\ge \frac{\varepsilon^4}{16}|V_t|$, such that $|d(X,Y) - d(V_s,V_t)| \ge
\varepsilon^4$.
\item {\bf Count regular pairs:} If there are at most $\varepsilon k_i^2$ pairs that are not verified as $\varepsilon$-regular, then halt. $P_i$ is an $\varepsilon$-regular partition.
\item {\bf Refinement:} Otherwise apply the Refinement Algorithm and
Lemma \ref{lnew2}, where $P = P_i, k = k_i, \gamma =
\frac{\varepsilon^4}{16}$, and obtain a partition $Q$ with $1 +
k_i4^{k_i}$ classes.
\item {\bf Iteration:} Let $k_{i+1} = k_i4^k_i, P_{i+1} = Q, i = i+1$, and go to step 2.
\end{enumerate}
}

Since the index cannot exceed $1/2$, the algorithm must halt after
at most $\lceil 10\gamma^{-5} \rceil$ iterations (see \cite{AD}).
Unfortunately, in each iteration the number of classes increases to
$k4^k$ from $k$. This implies that the graph $G$ must be indeed
astronomically large (a tower function) to ensure the completion of
this procedure. As mentioned before, Gowers \cite{G} proved that
indeed this tower function is necessary in order to guarantee an
$\varepsilon$-regular partition for {\em all} graphs. The size requirement
of the algorithm above makes it impractical for real world
situations where the number of vertices typically is a few thousand.

\subsubsection{Frieze-Kannan Version}

The Frieze-Kannan constructive version is quite similar to the
above, the only difference is how to check regularity of the pairs
in Step 2. Instead of Lemma \ref{lnew1}, another lemma is used based
on the computation of singular values of matrices. For the sake of
completeness we present the details.

\begin{lemma}[Frieze-Kannan \cite{FK}]\pont\label{singularFK}
Let $W$ be an $R \times C$ matrix with $|R| = p$ and $|C| = q$ and
$W_\infty\leq 1$ and $\gamma$ be a positive real.

\begin{enumerate}
\item[a] If there exist $S \subseteq R, T \subseteq C $ such that $|S| \geq \gamma p, |T| \geq \gamma q$ and $|W(S,T)| \geq \gamma |S| |T|$
then $\sigma_1(W) \geq \gamma^3 \sqrt{pq}$ (where $\sigma_1$ is the
first singular value).
\item[b] If $\sigma_1(W) \geq \gamma \sqrt{pq}$ then there exist  $S \subseteq R, T \subseteq C $ such that $|S| \geq \gamma'p, |T| \geq \gamma'q$
and $W(S,T) \geq \gamma' |S| |T|$, where $\gamma' =
\frac{\gamma^3}{108}$. Furthermorem $S,T$ can be constructed in
polynomial time.
\end{enumerate}

\end{lemma}

Combining Lemmas \ref{lnew2} and \ref{singularFK}, we get an
algorithm for finding an $\varepsilon$-regular partition, quite
similar to the Alon {\em et al.} version \cite{AD}, which we present
below:

{\bf Regular Partitioning Algorithm (Frieze-Kannan):} {\em Given a
graph $G$ and $\eps$, construct a $\eps$-regular partition.
\begin{enumerate}
\item {\bf Initial partition:} Arbitrarily divide the vertices of $G$ into an equitable partition $P_1$ with classes $V_0, V_1, \ldots, V_b$, where
$|V_1| = \lfloor \frac{n}{b} \rfloor$ and hence $|V_0| < b$. Denote
$k_1 = b$.
\item {\bf Check regularity:} For every pair $(V_s,V_t)$ of $P_i$, compute $\sigma_1(W_{r,s})$. If the pair $(V_r, V_s)$ are not $\varepsilon$-regular
then by Lemma \ref{singularFK} we obtain a proof that they are not
not $\gamma = \varepsilon^9/108$-regular.
\item {\bf Count regular pairs:} If there are at most $\eps k_i^2$ pairs that produce proofs of non $\gamma$-regularity, then halt. $P_i$ is an $\eps$-regular partition.
\item {\bf Refinement:} Otherwise apply the Refinement Algorithm and
Lemma \ref{lnew2}, where $P = P_i, k = k_i, \gamma =
\frac{\eps^9}{108}$, and obtain a partition $P'$ with $1 +
k_i4^{k_i}$ classes.
\item {\bf Iteration:} Let $k_{i+1} = k_i4^{k_i}, P_{i+1} = P', i = i+1$, and go to step 2.
\end{enumerate}
}

This algorithm is guaranteed to finish in at most
$\varepsilon^{-45}$ steps with an $\varepsilon$-regular partition.

\section{Modifications to the Constructive Version}\label{mod}

We see that even the constructive versions are not directly
applicable to real world scenarios. We note that the above
algorithms have such restrictions because their aim is to be
applicable to {\em all} graphs. Thus, to make the regularity lemma
truly applicable we would have to give up our goal that the lemma
should work for {\em every} graph and should be content with the
fact that it works for {\em most} graphs. To ensure that this
happens, we modify the Regular Partitioning Algorithm(s) so that
instead of constructing a regular partition, we find an {\em
approximately} regular partition, which should be much easier to
construct. We have the following 3 major modifications to the
Regular Partitioning Algorithm.

{\bf Modification 1:} We want to decrease the cardinality of atoms
in each iteration. In the above Refinement Algorithm the cardinality
of the atoms may be $2^{k-1}$, where $k$ is the number of classes in
the current partition. This is because the algorithm tries to find
all the possible $\varepsilon$-irregular pairs such that this
information can then be embedded into the subsequent refinement
procedure. Hence potentially each class may be involved with up to
$(k-1)$ $\varepsilon$-irregular pairs. One way to avoid this problem
is to bound this number. To do so, instead of using all the
$\varepsilon$-irregular pairs, we only use some of them.
Specifically, in this paper, for each class we consider at most one
$\varepsilon$-irregular pair that involves the given class. By doing
this we reduce the number of atoms to at most $2$. We observe that
in spite of the crude approximation, this seems to work well in
practice.

{\bf Modification 2:} We want to bound the rate by which the class
size decreases in each iteration. As we have at most $2$ atoms for
each class, we could significantly increase $m$ used in the
Refinement Algorithm as $m = \frac{|V_i|}{l}$, where a typical value
of $l$ could be $3$ or $4$, much smaller than $4^k$. We call this
user defined parameter $l$ the refinement number.

{\bf Modification 3:} Modification 2 might cause the size of the
exceptional class to increase too fast. Indeed, by using a smaller
$l$, we risk putting $\frac{1}{l}$ portion of all vertices into
$V_0$ after each iteration. To overcome this drawback, we
``recycle'' most of $V_0$, i.e. we move back most of the vertices
from $V_0$. Here is the modified Refinement Algorithm.

{\bf Modified Refinement Algorithm:} {\em Given a $\gamma$-irregular
equitable partition $P$ of the vertex set $V=V_0\cup V_1\cup \ldots
\cup V_k$ with $\gamma = \frac{\varepsilon^4}{16}$ and refinement
number $l$, construct a new partition $Q$.\\
For each pair $(V_s, V_t)$, $1 \leq s < t \leq k$, we apply Lemma
\ref{lnew1} with $A=V_s$, $B=V_t$ and $\varepsilon$. For a fixed $s$
if $(V_s, V_t)$ is found to be $\varepsilon$-regular for all $t\not=
s$ we do nothing, i.e. $V_s$ is one atom. Otherwise, we select one
$\varepsilon$-irregular pair $(V_s, V_t)$ randomly and the
corresponding certificate partitions $V_s$ into two atoms. Set $m =
\lfloor \frac{|V_i|}{l}\rfloor$, $1 \le i \le k$. Then first we
choose a maximal collection $Q'$ of pairwise disjoint subsets of $V$
such that every member of $Q'$ has cardinality $m$ and every atom
$A$ contains exactly $\lfloor \frac{|A|}{m}\rfloor$ members of $Q'$.
Then we unite the leftover vertices in each $V_s$, if there are at
least $m$ vertices then we select one more subset of size $m$ from
these vertices, we add these sets to $Q'$ and finally we add all
remaining vertices to the exceptional class, resulting in the
partition $Q$. The collection $Q$ is an equitable partition of $V$
into at most $1+lk$ classes. }

\noindent Now, we present our modified Regular Partitioning
Algorithm. There are three main parameters to be selected by the
user: $\varepsilon$, the refinement number $l$ and $h$, the minimum
class size when we must halt the refinement procedure. The parameter
$h$ is used to ensure that if the class size has gone too small then
the procedure should not continue.

{\bf Modified Regular Partitioning Algorithm (or the Practical
Regularity Partitioning Algorithm):} {\em Given a graph $G$ and
parameters $\varepsilon$, $l$, $h$, construct an approx.
$\varepsilon$-regular partition.
\begin{enumerate}
\item {\bf Initial partition:} Arbitrarily divide the vertices of $G$ into an equitable partition $P_1$ with classes $V_0, V_1, \ldots, V_l$, where
$|V_1| = \lfloor \frac{n}{l} \rfloor$ and hence $|V_0| < l$. Denote
$k_1 = l$.
\item {\bf Check size and regularity:} If $|V_i| < h$, $1\leq i \leq k$, then halt.
Otherwise for every pair $(V_s,V_t)$ of $P_i$, verify if it is $\varepsilon$-regular or find
$X \subset V_s, Y \subset V_t, |X| \ge \frac{\varepsilon^4}{16}|V_s|,|Y|
\ge \frac{\varepsilon^4}{16}|V_t|$, such that $|d(X,Y) - d(V_s,V_t)| \ge
\varepsilon^4$.
\item {\bf Count regular pairs:} If there are at most $\varepsilon k_i^2$ pairs that are not verified as $\varepsilon$-regular, then halt. $P_i$ is an $\varepsilon$-regular partition.
\item {\bf Refinement:} Otherwise apply the Modified Refinement Algorithm,
where $P = P_i, k = k_i, \gamma = \frac{\varepsilon^4}{16}$, and obtain a partition $Q$ with $1 + lk_i$ classes.
\item {\bf Iteration:} Let $k_{i+1} = l k_i, P_{i+1} = Q, i = i+1$, and go to step 2.
\end{enumerate}
}

The Frieze-Kannan version is modified in an identical way.

\section{Application to Clustering}\label{app}

To make the regularity lemma applicable in clustering settings, we
adopt the following two phase strategy (as in \cite{SP}):
\begin{enumerate}
\item {\bf Application of the Practical Regularity Partitioning Algorithm:} In the first
stage we apply the Practical Regularity partitioning algorithm as
described in the previous section to obtain an approximately regular
partition of the graph representing the data. Once such a partition
has been obtained, the reduced graph as described in Definition
\ref{reduced} could be constructed from the partition.

\item {\bf Clustering the Reduced Graph:} The reduced graph as constructed above
would preserve most of the properties of the original graph (see
\cite{KSSS}). This implies that any changes made in the reduced
graph would also reflect in the original graph. Thus, clustering the
reduced graph would also yield a clustering of the original graph.
We apply spectral clustering (though any other pairwise clustering
technique could be used, e.g. in \cite{SP} the dominant-set
algorithm is used) on the reduced graph to get a partitioning and
then project it back to the higher dimension. Recall that vertices
in the exceptional set $V_0$ are leftovers from the refinement
process and must be assigned to the clusters obtained. Thus in the
end these leftover vertices are redistributed amongst the clusters
using a k-nearest neighbor classifier to get the final grouping.

\end{enumerate}

\section{Empirical Validation}\label{test}
In this section we present extensive experimental results to
indicate the efficacy of this approach by employing it for
clustering on a number of benchmark datasets. We also compare the
results with spectral clustering in terms of accuracy. We also
report results that indicate the amount of compression obtained by
constructing the reduced graph. As discussed later, the results also
directly point to a number of promising directions of future work.
We first review the datasets considered and the metrics used for
comparisons.

\subsection{Datasets and Metrics Used} \label{datasets}
The datasets considered for empirical validation were taken from the
University of California, Irvine machine learning repository
\cite{UCI}. A total of 12 datasets were used for validation. We
considered datasets with real valued features and associated labels
or ground truth. In some datasets (as described below) that had a
large number of real valued features, we removed categorical
features to make it easier to cluster. Unless otherwise mentioned,
the number of clusters was chosen so as to equal to the number of
classes in the dataset (i.e. if the number of classes in the ground
truth is 4, then the clustering results are for k = 4). An attempt
was made to pick a wide variety of datasets, i.e. with integer
features, binary features, synthetic datasets and of course real
world datasets with both very high and small dimensionality.

The following datasets were considered (for details about the
datasets see \cite{UCI}): (1) Red Wine (R-Wine) and (2) White Wine
(W-Wine), (3) The Arcene dataset (Arcene), (4) The Blood Transfusion
Dataset (Blood-T), (5) The Ionosphere dataset (Ionos), (6) The
Wisconsin Breast cancer dataset (Cancer), (7) The Pima Indian
diabetes dataset (Pima), (8) The Vertebral Column dataset
(Vertebral-1), the second task (9) (Vertebral-2) is considered as
another dataset, (10) The Steel Plates Faults Dataset (Steel), (11)
The Musk 2 (Musk) dataset and (12) Haberman's Survival (Haberman)
data.

Next we discuss the metric used for comparison with other clustering
algorithms. For evaluating the quality of clustering, we follow the
approach of \cite{WuScho} and use the cluster accuracy as a measure.
The measure is defined as: $$\displaystyle Accuracy = 100*
\biggl(\frac{\sum_{i=1}^{n} \delta (y_i, map(c_i))}{n} \biggr ), $$
where $n$ is the number of data-points considered, $y_i$ represents
the true label (ground truth) while $c_i$ is the obtained cluster
label of data-point $x_i$. The function $\delta(y,c)$ equals 1 if
the true and the obtained labels match ($y=c$) and 0 if they don't.
The function $map$ is basically a permutation function that maps
each cluster label to the true label. An optimal match can be found
by using the Hungarian Method for the assignment problem
\cite{Kuhn}.

\subsection{Case Study} \label{casestudy}
Before reporting comparative results on benchmark datasets, we first
consider one dataset as a case study. While experiments reported in
this case study were carried on all the benchmark datasets
considered, the purpose here is to illustrate the investigations
conducted at each stage of application of the regularity lemma. An
auxiliary purpose is also to underline a set of guidelines on what
changes to the practical regularity partitioning algorithm proved to
be useful.

For this task we consider the Red Wine dataset which has 1599
instances with 11 attributes each, the number of classes involved is
six. It must be noted though that the class distribution in this
dataset is pretty skewed (with the various classes having 10, 53,
681, 638, 199 and 18 datapoints respectively), this makes clustering
this dataset quite difficult when k = 6. We however consider both k
= 6 and k = 3 to compare results with spectral clustering.

\begin{figure}
\centering
\includegraphics[width=2in]{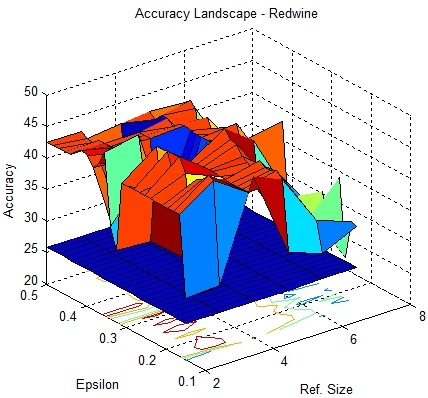}%
\hspace{0.5in}%
\includegraphics[width=2in]{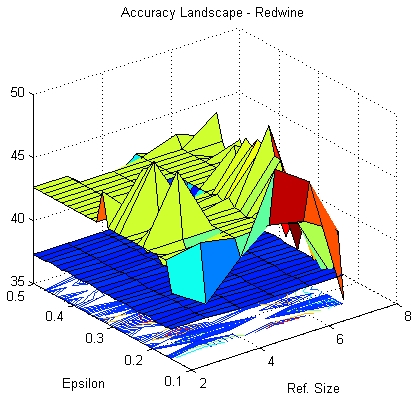}
\caption{Accuracy Landscape for Regularity Clustering on the Red Wine Dataset for different values of $\varepsilon$ and refinement size $l$ (with k = 6 on the left and k = 3 on the right). The Plane cutting through in blue represents accuracy by running self-tuned spectral clustering using the fully connected similarity graph.}
\label{fig:CaseStudy}
\end{figure}

Recall that our method has two meta-parameters that need to be user
specified (or estimated by cross-validation): $\varepsilon$ and $l$.
Note that $h$ is usually decided so that it is at least as big as
$\frac{1}{\varepsilon}$. The first set of experiments thus explore
the accuracy landscape of regularity clustering spanned over these
two parameters. We consider 25 linearly spaced values of
$\varepsilon$ between 0.15 and 0.50. The refinement number $l$, as
noted in Section \ref{mod}, can not be too large. Since it can only
take integer values, we consider six values from 2 to 7. For the
sake of comparison, we also obtain clustering results on the same
dataset with spectral clustering with self tuning \cite{selftuned}
(both using all connected and k-nearest neighbour graph versions)
and k-means clustering. Figure \ref{fig:CaseStudy} gives the
accuracy of the Regularity Clustering on a grid of $\varepsilon$ and
$l$. Even though this plot is only for exploratory purposes, it
shows that the accuracy landscape is in general much better than the
accuracy obtained by spectral clustering for this dataset.


An important aspect of the Regularity Clustering method is that by
using a modified constructive version of the Regularity Lemma we
obtain a much reduced representation of the original data. The size
of the reduced graph depends both on $\varepsilon$ and $l$. However,
in our observation it is more sensitive to changes to $l$ and
understandably so. From the grid for $\varepsilon$ and $l$ we take
three rows to illustrate the obtained sizes of the reduced graph
(more precisely, the dimensions of the affinity matrix of the
reduced graph). We compare these numbers with the original dataset
size. As we note in the results over the benchmark datasets in
section \ref{benchmark}, this compression is quite big in larger
datasets.

\begin{table}
\caption{Reduced Graph Sizes. Original Affinity Matrix size : 1599 $\times$ 1599}
\centering

\begin{tabular}{| c | c | c | c | c | c | c |}
\hline \hline
\backslashbox{{\bf $\varepsilon$}}{{\bf $l$}} & {\bf 2} & {\bf 3} & {\bf 4} & {\bf 5} & {\bf 6} & {\bf 7} \\ \hline
{\bf 0.15} & 16 $\times$ 16 & 27 $\times$ 27 & 27 $\times$ 27 & 27 $\times$ 27  & 36 $\times$ 36 & 49 $\times$ 49 \\ \hline
{\bf 0.33} & 49 $\times$ 49 & 49 $\times$ 49 & 66 $\times$ 66 & 66 $\times$ 66  & 66 $\times$ 66 & 66 $\times$ 66 \\ \hline
{\bf 0.50} & 66 $\times$ 66 & 66 $\times$ 66  & 66 $\times$ 66 & 66 $\times$ 66 & 66 $\times$ 66 & 66 $\times$ 66 \\ \hline
\end{tabular}
\label{redgraphsizes}
\end{table}

The proof of the Regularity Lemma is using a potential function, the
index of the partition defined earlier in Definition
\ref{potential}. In each refinement step the index increases
significantly. Surprisingly this remains true in our modified
refinement algorithm when the number of partition classes is not
increasing as fast as in the original version, see Table
\ref{potentialincrease}. Another interesting observation is that if
we take $\varepsilon$ sufficiently high, we do get a regular
partition in just a few iterations. A few examples where this was
noticed in the Red Wine dataset are mentioned in Table
\ref{regpartitions}.

\begin{table}
\caption{Illustration of Increase in Potential}
\centering
\begin{tabular}{| c | c | c | c | c |}
\hline \hline
\backslashbox{{ \bf ($\varepsilon$,$l$)}}{{\bf $ind(P)$}} & {\bf $ind(P_1)$} & {\bf $ind(P_2)$} & {\bf $ind(P_3)$} & {\bf $ind(P_4)$} \\ \hline
{\bf 0.15, 2} & 0.1966 & 0.2892 & 0.3321 & 0.3539  \\ \hline
{\bf 0.33, 2} & 0.1966 & 0.2883 & 0.3321 & 0.3683   \\ \hline
{\bf 0.50, 2} & 0.1965 & 0.2968 & 0.3411 & 0.3657    \\ \hline
\end{tabular}
\label{potentialincrease}
\end{table}

\begin{table}
\caption{Regular Partitions with req. no. of regular pairs and actual no. present}
\centering
\begin{tabular}{| c | c | c | c |}
\hline \hline
($\varepsilon$, $l$) & \# for $\varepsilon$-regularity & \# of Reg. Pairs & \# Iterations \\ \hline
{\bf 0.6, 2} & 1180 & 1293 & 6 \\ \hline
{\bf 0.7, 6} & 352 & 391 & 2 \\ \hline
{\bf 0.7, 7} & 506 & 671  & 2 \\ \hline
\end{tabular}
\label{regpartitions}
\end{table}

Finally, before reporting results we must make a comment on
constructing the reduced graph. The reduced graph was defined in
Definition \ref{reduced}. But note that there is some ambiguity in
our case when it comes to constructing the reduced graph. The
reduced graph $G^R$ is constructed such that the vertices correspond
to the classes in the partition and the edges are associated to the
$\varepsilon$-regular pairs between classes with density above $d$.
However, in many cases the number of regular pairs is quite small
(esp. when $\varepsilon$ is small) making the matrix too sparse,
making it difficult to find the eigenvectors. Thus for technical
reasons we added all pairs to the reduced graph. We contend that
this approach works well because the classes that we consider (and
thus the densities between them) are obtained after the modified
refinement procedure and thus enough information is already embedded
in the reduced graph.

\subsection{Clustering Results on Benchmark Datasets}\label{benchmark}

In this section we report results on a number of datasets described
earlier in Section \ref{datasets}. We do a five fold
cross-validation on each of the datasets, where a validation set is
used to learn the meta-parameters for the data. The accuracy
reported is the average clustering quality on the rest of the data
after using the learned parameters from the validation set. We use a
grid-search to learn the meta-parameters. Initially a coarse grid is
initialized with a set of 25 linearly spaced values for
$\varepsilon$ between 0.15 and 0.50 (we don't want $\varepsilon$ to
be outside this range). For $l$ we simply pick values from 2 to 7
simply because that is the only practical range that we are looking
at.

\begin{table}
\caption{Clustering Results on UCI Datasets. Regular1 and Regular2
represent the results by the versions due to Alon {\em et al.} and
Frieze-Kannan, respectively. Spect1 and Spect2 give the results for
spectral clustering with a k-nearest neighbour graph and a fully
connected graph, respectively. Follow the text for more details.}
\centering
\begin{tabular}{llllllll}
\hline\noalign{\smallskip}
{\bf Dataset} & {\bf \# Feat.} & {\bf Comp.} &  {\bf Regular1} &  {\bf Regular2} & {\bf Spect1} & {\bf Spect2} & {\bf k-means}\\
\noalign{\smallskip}
\hline
\noalign{\smallskip}
R-Wine & 11 & 1599-49 & {\bf 47.0919} & {\bf 46.8342} & 23.9525 & 23.9524 & 23.8899  \\
W-Wine & 11 & 4898-125 &{\bf 44.7509} &{\bf 44.9121} & 23.1319  & 20.5798  & 23.8465   \\
Arcene & 10000 & 200-9 & {\bf 68} & {\bf 68} & 61  & 62  & 59  \\
Blood-T & 4 & 748-49 & {\bf 76.2032} & {\bf 75.1453} & 65.1070  & 66.2331  & 72.3262   \\
Ionos & 34 & 351-25 & {\bf 74.0741} & {\bf 74.6787} & 70.0855  & 70.6553  & 71.2251   \\
Cancer & 9 & 683-52 & 93.5578 & 93.5578 & {\bf 97.2182}   & 97.2173   & 96.0469   \\
Pima & 8 & 768-52 & {\bf 65.1042 } & {\bf 64.9691} &  51.5625  & 60.8073 & 63.0156   \\
Vertebral-1 & 6 & 310-25 &67.7419 & 67.8030 & {\bf 74.5161} & 71.9355 & 67.0968    \\
Vertebral-2 & 6 & 310-25 &{\bf 70 } &{\bf 69.9677} & 49.3948 & 48.3871  & 65.4839   \\
Steel & 27 & 1941-54 & {\bf 42.5554} & {\bf 43.0006} & 29.0057   & 34.7244  &  29.7785  \\
Musk & 166 & 6598-126 & {\bf 84.5862} & {\bf 81.4344} & 53.9103   & 53.6072  &  53.9861  \\
Haberman & 3 & 306-16 & {\bf 73.5294} & {\bf 70.6899} & 52.2876   & 51.9608  &  52.2876  \\
\hline
\end{tabular}\label{UCIResultTable}
\end{table}

We compare our results with a fixed $\sigma$ spectral clustering
with both a fully connected graph (Spect2) and a k-nearest neighbour
graph (Spect1). For the sake of comparison we also include results
for k-means on the entire dataset. We also report results on the
compression that was achieved on each dataset in Table
\ref{UCIResultTable} (The compression is indicated in the format x-y
where x represents one dimension of the adjacency matrix of the
dataset and y of the reduced graph).

In the results we observe that the Regularity Clustering method, as
indicated by the clustering accuracies is quite powerful; it gave
significantly better results in 10 of the 12 datasets. It was also
observed that the regularity clustering method did not appear to
work very well in synthetic datasets. This seems understandable
given the quasi-random aspect of the Regularity method. We also
report that the results obtained by the Alon {\em et al.} and by the
Frieze-Kannan versions are virtually identical, which is not
surprising.

\section{Future Directions} \label{future}

We believe that this work opens up a lot of potential research
problems. First and foremost would be establishing theoretical
results for quantifying the approximation obtained by our
modifications to the Regularity Lemma. Also, the original Regularity
Lemma is applicable only while working with dense graphs. However,
there are sparse versions of the Regularity Lemma. These sparse
versions could be used in the first phase of our method such that
even sparse graphs (k-nearest neighbor graphs) could be used for
clustering, thus enhancing its practical utility even further.

A natural generalization of pairwise clustering methods leads to
hypergraph partitioning problems \cite{Bulo}, \cite{Zhou}. There are
a number of results that extend the Regularity Lemma to hypergraphs
\cite{ChungHyper}, \cite{GowersHyper}, \cite{RodlHyper2}. It is thus
natural that our methodology could be extended to hypergraphs and
then used for hypergraph clustering.

In final summary, our work gives a way to harness the Regularity
Lemma for the task of clustering. We report results on a number of
benchmark datasets which strongly indicate that the method is quite
powerful. Based on this work we also suggest a number of possible
avenues for future work towards improving and generalizing this
methodology.


\begin{thebibliography}{4}
\bibitem{AD} N. Alon, R.A. Duke, H. Lefmann, V. R\"{o}dl, R. Yuster, The Algorithmic Aspects of the Regularity Lemma.
In: J. Alg., 16, pp. 80-109, (1994).

\bibitem{Bulo} S. Bulo, M. Pelillo, A Game-Theoretic Approach to Hypergraph Clustering. In: NIPS,
22, pp. 1571-1579, (2009).

\bibitem{ChungHyper} F. Chung, Regularity Lemmas for Hypergraphs and Quasi-randomness. In: Random Struct.
Alg., 2, pp. 241-252, (1991).

\bibitem{ET} P. Erd\"os, P. Tur\'an, On Some Sequences of Integers. In: J. London Math. Soc, 11, pp. 261-264,
(1936).

\bibitem{FowlChung} C. Fowlkes, S. Belongie, F. Chung, J. Malik, Spectral grouping using the Nystr\"{o}m Method. In: IEEE Trans. PAMI, 26, pp. 214-225,
(2004).

\bibitem{FK} A. Frieze, R. Kannan, A Simple Algorithm for Constructing Szemer\'edi's Regularity Partition. In: Electron. J. Comb, 6 (1), R17, (1999).

\bibitem{UCI} A. Frank, A. Asuncion, UCI Machine Learning Repository, Irvine, CA: University of California, School of Information and Computer Science, (2010).

\bibitem{G} W.T. Gowers, Lower bounds of tower type for Szemer\'edi's uniformly lemma. In: Geom. Funct. Anal, 7, pp. 322-337,
(1997).

\bibitem{GowersHyper} W.T. Gowers, Hypergraph regularity and the Multidimensional Szemer\'edi theorem, In: Annals of Mathematics, (2) 166 no. 3, pp. 897-946,
(2007).

\bibitem{Abel} W.T. Gowers, The Work of Endre Szemer\'{e}di, Online at \url {http://www.abelprize.no/c54147/binfil/download.php?tid=54060}

\bibitem{GRSS1} A. Gy\'arf\'as, M. Ruszink\'o, G.N. S\'ark\"ozy, E. Szemer\'edi, Three-color Ramsey numbers for paths, In: Combinatorica, 27(1), pp. 35-69,
(2007).

\bibitem{KRT} Y. Kohayakawa, V. R\"{o}dl, L. Thoma, An optimal algorithm for checking regularity.In: SIAM J. Comput, 32(5), pp. 1210-1235,
(2003).

\bibitem{KSS5} J. Koml\'os, G.N. S\'ark\"ozy, E. Szemer\'edi, Blow-up Lemma, In: Combinatorica, 17(1), pp. 109-123,
(1997).

\bibitem{KSSS} J. Koml\'os, A. Shokoufandeh, M. Simonovits, E. Szemer\'edi, The Regularity Lemma and Its Applications in Graph Theory. In: Theoretical Aspects of Comp. Sci., LNCS 2292,
pp. 84-112, (2002).

\bibitem{Kuhn} H.W. Kuhn, The Hungarian method for the Assignment Problem, In: Naval Research Logistics, 52(1), 2005.
Originally appeared in Naval Research Logistics Quarterly, 2, 1955, pp. 83-97.

\bibitem{NJW} A. Ng, M. Jordan, Y. Weiss, On Spectral Clustering: Analysis and an algorithm. In T. Dietterich, S. Becker, and
Z. Ghahramani (Eds.),NIPS, MIT Press, 14, pp. 849-856, (2002).

\bibitem{RodlHyper2} V. R\"odl, B. Nagle, J. Skokan, M. Schacht, Y. Kohayakawa, The Hypergraph Regularity Method and its applications, In: PNAS, 102,
pp. 8109-8113, (2005).

\bibitem{shimalik} S. Shi, J. Malik, Normalized Cuts and Image Segmentation, In: IEEE Trans. PAMI, vol 22, no. 8, pp. 888-905,
(2000).

\bibitem{SP} A. Sperotto, M. Pelillo, Szemer\'{e}di's Regularity Lemma and its Applications to Pairwise Clustering and
Segmentation. In: EMMCVPR, LNCS 4679. Springer, (2007).

\bibitem{Sz} E. Szemer\'edi, Regular Partitions of Graphs, Colloques Internationaux C.N.R.S. $\mbox{N}^{\underline{o}}$ 260 -
Probl\`emes Combinatoires et Th\'eorie des Graphes, Orsay, pp.
399-401, (1976).

\bibitem{WuScho} M. Wu, B. Sch\"olkopf, A Local Learning Approach for Clustering. In: NIPS, pp. 1529-1536,
(2007).

\bibitem{selftuned} L. Zelnik-Manor, P. Perona, Self-tuning Spectral Clustering. In L. K. Saul, Y. Weiss, and L. Bottou, eds., NIPS, MIT
Press, Cambridge, MA,  pp. 1601-1608, (2005).


\bibitem{Zhou} D. Zhou, J. Huang, B. Sch\"olkopf, Learning with Hypergraphs: Clustering, Classification, and Embedding. In: NIPS, 19,
pp. 1601–1608, (2007).

\end{thebibliography}
\end{document}